\documentclass{amsart}
\usepackage{amssymb}
\usepackage{amsfonts}
\usepackage{amsmath}
\usepackage[legalpaper,bookmarks=true,colorlinks=true,linkcolor=blue,citecolor=blue]{hyperref}
\usepackage{graphicx}%
\usepackage{fancyhdr}
\usepackage{color}
\usepackage[mathlines]{lineno}
\usepackage{lscape}
\usepackage{epsfig}
\usepackage{natbib}
\usepackage{geometry}
\usepackage{tgbonum}
\usepackage[]{listings}
\fontfamily{qcr}\selectfont

\newtheorem{theorem}{Theorem}
\theoremstyle{plain}

\numberwithin{equation}{section}

\newcommand{\Bin}{\bigskip \noindent}

\newcommand{\Ni}{\noindent}

\begin{document}
\Large
\title[random analysis on the a generalized hyperbolic distribution L\'evy process]{Elements of random analysis about the gamma generalized hyperbolic distribution Levy stochastic process}

\author{Nafy Ngom}
\author{Aladji Babacar Niang}
\author{Soumaila Dembel\'e}
\author{Gane Samb Lo}

\begin{abstract} In this paper, we study some aspects on random analysis on the L\'evy stochastic processes with margins following generalized  hyperbolic distributions generated by gamma laws. In particular we study the boundedness of its total variations and the quadratic variations. Next we give an empirical construction that enables the graphical representation of the paths of such stochastic processes. Comparisons with the Brownian motions are considered.\\  

\noindent Nafy Ngom\\
LERSTAD, Gaston Berger University, Saint-Louis, S\'en\'egal.\\

\noindent Dr Aladji Babacar Niang\\
LERSTAD, Gaston Berger University, Saint-Louis, Senegal\\
HEC Faculty, UNIL, University of Lausanne, Suisse \\
Emails: niang.aladji-babacar@ugb.edu.sn, aladjibacar93@gmail.com, aladjibabacar.niang@unil.ch\\
Suisse address: 1027, Chemin de la Fauvette 2, Lonay\\
Senegal address: 088, HLM Grand Yoff, Dakar\\

\noindent Dr Soumaila Dembel\'e\\
Universit\'e des Sciences Sociale et de Gestion de Bamako (USSGB)\\
Facult\'e des Sciences \'Economiques et de Gestion (FSEG)\\
Email: dembele.soumaila@ugb.edu.sn, soumailadembeleussgb@gmail.com\\

\noindent $^{\dag}$ Gane Samb Lo.\\
LERSTAD, Gaston Berger University, Saint-Louis, S\'en\'egal (main affiliation).\newline
LSTA, Pierre and Marie Curie University, Paris VI, France.\newline
AUST - African University of Sciences and Technology, Abuja, Nigeria\\
gane-samb.lo@edu.ugb.sn, gslo@aust.edu.ng, ganesamblo@ganesamblo.net\\
Permanent address : 1178 Evanston Dr NW T3P 0J9,Calgary, Alberta, Canada.\\

\end{abstract}
\maketitle

\section{Introduction}

\Ni In the recent years, the generalized hyperbolic distribution (\textit{GHD}) is attracting much attention in various branches of Statistics and Probability theory, especially in statistical estimation and in finance modeling (see \cite{mcneil}). This class of laws is seen as a credible alternative to the Gaussian assumption of data both in skewness and in shape. Also, as infinitely divisible laws, the L\'evy processes they generate began to be studied (see \cite{applebaum}).\\

\Bin To be able to expose the aim of this paper and to describe its motivation, we need to develop a little the three aspects  aforementioned in the next line.\\

\subsection{Infinitely divisibility laws or random variables} \label{ss_01_01}$ $\\

\Ni The probability law of a real-valued random variable $X$ associated to the characteristic function $\psi_X$ is infinitely divisible, denoted as [\textit{idecomp}], if and only if it satisfies one of the criteria below.\\

\Ni \label{crit1} (Criterion 1) A random variable $X$ is infinitely divisible or infinitely decomposable (\textit{idecomp}) if and only if for any $n\geq 1$, it can be equal in law to a sum of $n$ independent and identically distributed random variables.\\

\Ni \label{crit2} (Criterion 2) A random variable $X$ of characteristic function $\psi_X$ (or of third moment generating function $\phi_X$ if $X\geq 0$)  is  \textit{idecomp} if and only if: for any $n\geq 1$, there exists a characteristic function $\psi$ (or a third moment generating function $\phi$ of a non-negative random variable) such that

$$
\psi_X=\psi^n \ \ (or \ \ \phi_X=\phi^n).
$$

\Bin \label{crit3} (Criterion 3) A random variable $X$ of characteristic function $\psi_X$ (or of third moment generating function $\phi_X$ if $X\geq 0$)  is  \textit{idecomp} if and only if: for any $n\geq 1$, $\psi_X^{1/n}$ is a characteristic function (or $\phi_X^{1/n}$ is a third moment generating function).\\

\Bin \textbf{Examples}. Let us consider three classical \textit{idecomp} probability laws: normal law $U_1= \mathcal{N}(m,\sigma^2)$, translation of a Poisson law $\mathcal{P}(\lambda)$ by $a\geq 0$: $U_2=a +\mathcal{P}(\lambda)=:\mathcal{P}(a,\lambda)$, and translation of a gamma law $\gamma(a,\beta)$ by $a\geq 0$: $U_3=a+\gamma(a,\beta)=:\gamma(a,a,\beta)$ with characteristic functions, for $t\in \mathbb{R}$,

\begin{eqnarray*}
&&\psi_{U_1}(t)=\exp(imt - \sigma^2t^2 /2)  , \ \ \ \psi_{U_2}(t)=\exp\left(iat + \lambda(e^{i t}-1)\right)\\
&& \psi_{U_3}(t)=e^{iat} \left(1-it/\beta\right)^{-a}.
\end{eqnarray*}

\Bin By taking the powers at $1/n$, $n\geq 1$, we get for any $t \in \mathbb{R}$,

$$
\psi_{U_1}(t)^{1/n}=\exp\left(i\{m/n\}t - (\sigma/\sqrt{n})^2t^2/2\right), \ \ \psi_{U_2}(t)^{1/n}=\exp\left(i\{a/n\}t + \{\lambda/n\}(e^{it}-1)\right),   
$$

$$
\ \ \  \ \ \psi_{U_3}(t)^{1/n}=e^{i\{a/n\}t} \left(1-it/\beta\right)^{\{-a/n\}}.
$$

\Bin By identifying the characteristic functions, we can see that the $\psi_{U_i}(t)^{1/n}$, $1\leq i \leq 3$, are \textit{cha.f}'s of random variables of law $\mathcal{N}(m/n,\sigma^2/n)$, $\mathcal{P}(a/n,\lambda/n)$ and $\gamma(a/n,a/n,\beta)$ respectively. So, these probability laws are \textit{idecomp}.\\

\Ni Full exposure of \textit{idecomp} laws are treated in \cite{loeve}, \cite{feller2}, \cite{zolotarev}, etc., and the readers are directed to them for getting deeper properties.\\ 

\subsection{Generalized Hyperbolic Distributions (\textbf{GHD})} \label{ss_01_02}$ $\\

\noindent A real valued-random variable $X$ is said to follow a normal variance mixture (\textit{NVM}) law if it is of the form

\begin{equation}
X=\mu +\sigma Z\sqrt{W}, \label{ghd-rep}
\end{equation}

 \bigskip \noindent where $\mu \in \mathbb{R}$, $\sigma>0$, $Z  \sim \mathcal{N}( 0,1)$, $Z$ independent of $W$, $W$ is some positive random variable.\newline

 \noindent As in \cite{mcneil}, we consider a few number of distributions for $W$:\newline

 \noindent a) The gamma distribution $\gamma(a, \beta)$ of parameters $a > 0$ and $\beta > 0$ of probability density function (\textit{pdf})

\begin{equation*}
f_{1}\left( x\right) =\frac{\beta ^{a }}{\Gamma \left( a \right) }
x^{a -1}e^{-\beta x}, \ x> 0;
\end{equation*}

\bigskip \noindent b) Inverse gamma distribution $Ig(a,\beta)$ of $a>0$ and $\beta>0$ of \textit{pdf} 

\begin{equation*}
f_{2}\left( x\right) =\frac{\beta ^{a }}{\Gamma \left( a \right) }
x^{-(a+1)} e^{-\beta/x}, \ \ x>0;
\end{equation*}

\bigskip \noindent c) The generalized inverse gamma $\textit{GIG}(a,b,c)$ distribution of parameters $a\geq 0$, $b\geq 0$, $c\geq 0$ of \textit{pdf}

\begin{equation*}
f_{3}(x) =\frac{1}{C\left( a,b,c\right) }x^{a-1}e^{-\left( bx+c/x\right) /2}, \ x>0,
\end{equation*}

\Bin where
 
\begin{equation*}
C\left( a,b,c\right) =\int_{0}^{+\infty}x^{a-1}e^{-\left( \left( bx+c/x\right)/2\right)} dx.
\end{equation*}

\Bin Some authors use the Bessel functions $K_a$ to express the constant of the distribution, and we have
 
\begin{equation*}
\frac{1}{C\left( a,b,c\right) }=\frac{b^{-a}\left( bc\right) ^{\frac{a}{2}}}{2K_{a}\left( \left( bc\right)^{\frac{1}{2}}\right) },
\end{equation*}

\bigskip \noindent with

\begin{equation*}
K_{a}\left( \left( bc\right) ^{\frac{1}{2}}\right) =\frac{c^{-a}\left(bc\right) ^{\frac{a}{b}}}{2}\int_{0}^{+\infty}x^{a-1}e^{-\left(
\left(bx+c/x\right) /2\right)} dx.
\end{equation*}

\bigskip \noindent The parameters $a$, $b$ and $c$ satisfy the conditions

\begin{eqnarray*}
a &=&0:b> 0\ \ and \ c> 0, \\
a &>&0:b> 0\ \ and \ c\geq 0, \\
a &< &0:b\geq 0 \ and \ c> 0.
\end{eqnarray*}

\bigskip \noindent We study the cases ($a> 0,b> 0,c=0$) and ($a< 0,b=0,c>0$) apart. The first case reduces to a gamma law $\gamma \left( a,b/2\right)$ of \textit{pdf}

\begin{equation*}
f_{4}\left( x\right) =\frac{(b/2)^{a}}{\Gamma \left( a\right)} x^{a-1}e^{-bx/2}, \ \ x\geq 0
\end{equation*}

\bigskip \noindent and the second to the following inverse gamma $Ig(a,c/2)$ distribution of \textit{pdf}

\begin{equation*}
f_{2,1}\left( x\right) =\frac{1}{D\left( a,c\right) }x^{-\left( a+1\right) }e^{-c/2x}, \ \ x> 0,
\end{equation*}

\bigskip \noindent with $a>0$, $c> 0$, and

$$
D\left(a,c\right)=\int_{0}^{+\infty} x^{-\left(a+1\right)} e^{-c/2x} \ dx.
$$ 

\bigskip \noindent We also have

\begin{equation*}
D\left( a,c\right) =\left( \frac{c}{2}\right) ^{-a}\Gamma ^{-}\left(
a\right) ,\ \ \Gamma ^{-}\left( a\right) =\int_{0}^{+\infty }x^{-\left(
a+1\right) }e^{-c/2x} \ dx,\ \ a>0, c>0.
\end{equation*}

\Bin If we take as $W$ one of this sample of distributions, the \textit{NVM} law becomes the Generalized Hyperbolic Distribution (\textit{GHD}). The \textit{GHD} model is full if we have $a\neq 0$, $b>0$ and $c>0$.\\

\Ni For any random variable $W>0$ with the density $f_{W}$, the \textit{pdf} of a \textit{NVM} random variable can be obtained by the classical diffeomorphical change of variables (see \cite{ips-mfpt-ang}, Chapter 2, Section 7, Part II) as follows. Let us consider the following diffeomorphism

\begin{eqnarray*}
h &:&\mathbb{R}\times \mathbb{R}_{+}^{\ast }\rightarrow \mathbb{R}\times 
\mathbb{R}_{+}^{\ast } \\
\left( z,w\right) &\longmapsto &\left( zw^{1/2},w\right)
=\left(u,v\right).
\end{eqnarray*}

\Bin The new \textit{pdf} of $\left( U,V\right) =h\left(Z, W\right)$ is given by

\begin{equation*}
f_{U,V}\left( u,v\right) =f_{Z}\left( uv^{-1/2}\right)
f_{W}\left(v\right) v^{-1/2} 1_{\mathbb{R}\times \mathbb{R}_{+}^{\ast }}\left(u,v\right).
\end{equation*}

\Bin The \textit{pdf} of $ZW^{1/2}$ is then

\begin{equation*}
f_{ZW^{1/2}}\left( u\right) =\frac{1}{\sqrt{2\pi }}\int_{0}^{+\infty} v^{-1/2} \ \exp(-u^2/(2v)) \ f_{W}\left( v\right) dv,\ \ u\in \mathbb{R}.
\end{equation*}

\Bin We are going to apply that formula for the distributions given above.\\

\Ni The distributions of $ZW^{1/2}$ corresponding to the \textit{pdf} $f_i$, $i \in \{1,2,3\}$ of $W$ are given as follows. Here we do not give the computations related to the application of the diffeomorphic transform, in particular their jacobian:

\begin{eqnarray*}
f_{g,ZW^{1/2}}(u)=\frac{\beta^a}{\sqrt{2\pi}\Gamma(a)} C\left(a-1/2, 2\beta,u^2\right), \ u \in \mathbb{R},
\end{eqnarray*}

\begin{eqnarray*}
f_{ig,ZW^{1/2}}(u)=\frac{2^{a+1/2}\beta ^{a}\Gamma(a+1/2)}{\Gamma \left( a\right) \sqrt{2\pi }} 
\left(u^{-2}+2\beta\right)^{-(a+1/2)}, \ u \in \mathbb{R},
\end{eqnarray*}

\begin{eqnarray*}
f_{gig,ZW^{1/2}}(u)=\frac{1}{\sigma C\left( a,b,c\right) \sqrt{2\pi}} C\left(a-1/2,b,c+u^2\right), \ \ u \in \mathbb{R}.
\end{eqnarray*}

\Bin Hence we deduce the \textit{GHD} associated to each case given by: \\

\Bin \textbf{1.} For $W \sim \gamma(a,\beta)$, $a>0$, $b>0$, we have

\begin{equation*}
f_{g,X}\left( u\right) =\frac{\beta^a}{\sigma \sqrt{2\pi}\Gamma(a)} C\left(a-1/2, 2\beta, \sigma^{-2}(u-\mu)^2\right), \ u\in \mathbb{R};
\end{equation*}

\Bin \textbf{2.}  For $W \sim Ig(a,\beta)$, $a>0$, $\beta>0$, we have

\begin{equation*}
f_{ig,X}\left( u\right) =\frac{2^{a+1/2}\beta ^{a}\Gamma(a+1/2)}{\sigma \Gamma \left( a\right) \sqrt{2\pi }}\left(\sigma^{-2}(u-\mu)^2+2\beta\right)^{-(a+1/2)}, \ u\in \mathbb{R};
\end{equation*}

\Bin \textbf{3.} For $W \sim Gig(a,b,c)$, $a>0$, $b>0$, $c>0$, 

\begin{equation*}
f_{gig,X}\left( u\right) =\frac{1}{\sigma C\left( a,b,c\right) \sqrt{2\pi}} C\left(a-1/2,b,c+\sigma^{-2}(u-\mu)^2\right), \ u\in \mathbb{R}.
\end{equation*}

\Bin To find the characteristic function $\psi_X$ of $X$, we denote by $H_{W}(t)=\mathbb{E}\left(e^{-tW}\right)$, $t\geq 0$, the moment function of $W$.

\Bin In general, by using the mathematical expectation tool, we get

\begin{eqnarray*}
\psi_X(t)&=&\mathbb{E}\left(\exp(i\mu t + i\sigma t W^{1/2} Z)\right)\\
&=&\exp(i\mu t) \ \mathbb{E}\biggr(\mathbb{E} \biggr(\exp\left(i \sigma t W^{1/2} Z\right) \biggr| W\biggr)\biggr)\\
&=&\exp(i\mu t) \ \mathbb{E}\biggr(\exp(-\sigma^2 t^2 W/2)\biggr)\\
&=& \exp(i\mu t) \int_{0}^{+\infty} \exp(-\sigma^2 t^2 w/2) f_W(w) \ dw\\
&=& \exp(i\mu t) \ H_{W}(\sigma^2t^2/2).
\end{eqnarray*}

\Bin  But for our three choices of $W$, \textit{i.e.} $W \in \left\{\gamma(a,\beta), Ig(a,\beta), \ gig(a,b,c)\right\}$, we have, respectively

$$
H_{g}(t)=(1+t/\beta)^{-a},  \  H_{Ig}(t)=\frac{\beta^a}{\Gamma(a)} C(-a,2t,2\beta)   \ and \ H_{gig}(t)=\frac{C(a,b+2t,c)}{C(a,b,c)}, \ t> 0.
$$

\Bin So for $H_W \in \{H_{g}, H_{Ig}, H_{gig}\}$, we have

\begin{eqnarray*}
\psi_X(t)&=&\exp(i\mu t) \ H_{W}(\sigma^2t^2/2), \ t>0.
\end{eqnarray*}

\Bin This may help showing that  $X$ is infinitely decomposable (\textit{idecomp}). Indeed, for any $k\geq 1$, we have

\begin{eqnarray*}
\psi_X(t)^{(1/k)}&=&\exp(i(\mu/k) t) \ H_{W}^{(1/k)}(\sigma^2t^2/2), \ t>0,
\end{eqnarray*}

\Bin and $X$ is \textit{idecomp} whenever $W$ is \textit{idecomp}. But it is possible, as in \cite{grosswald}, to show that any \textit{GHD} is \textit{decomp}, but this task is out of the scope of the paper. Since we focus on the \textit{GHD} associated to the gamma law, the characteristic function of $X$ is

\begin{equation}
\forall t \in \mathbb{R}, \ \  \psi_X(t)=\exp(i\mu t) \biggr(1 +\frac{\sigma^2t^2}{2\beta}\bigg)^{-a}, \label{gammagh-fchaf}
\end{equation}

\Bin and the application of Criterion 3, in page \pageref{crit3}, ensures that $Z$ is  \textit{idecomp}. \\

\Ni Generally,  for $W$ following the $\gamma(a,b)$ law, the $Ig(a,c)$ law or the $gig(a,b,c)$ law, $a>0$, $b>0,c>0$, we denote the \textit{GHD} respectively as the \textit{gamma-gh}$(a,b,0,\mu,\sigma)$ law, or the \textit{ig-gh}$(a,0,c,\mu,\sigma)$ law or the \textit{gig-gh}$(a,b,c,\mu,\sigma)$ law.\\

\Ni In this paper, we only treat with, \textit{gamma-gh}$(a,b,0,\mu,\sigma)$,  the \textit{gamma DH} distribution of parameters $a>0$,
$b>0$, $c=0$, $\sigma>0$ and $\mu \in \mathbb{R}$.\\

\Ni Our main source on \textit{GHD} distribution is \cite{mcneil}.\\

\subsection{L\'evy processes} \label{ss_01_03}$ $\\

\Ni By definition, a L\'evy process $\biggr(\Omega,\mathcal{A}, \mathbb{P}, (K_t)_{(t\geq 0)}\biggr)$ is defined as a stochastic process  with marginal laws $\{K(t), \ t\geq 0\}$ such that: \\

\Ni (ML) $\forall t\geq 0$, \ $L(t)\sim K(t)$, with $L(0)=0$ \textit{a.s}; \\

\Ni (II) $\forall k\geq 0$, \ $\forall \ 0=t_0<t_1<\cdots <t_k$, the \textit{r.v}'s $L(t_j)-L(t_{j-1})$ are independent; \\

\Ni (IS) $\forall 0\leq s<t$, \ $L(t)-L(s)\sim K(t-s)$; \\

\Ni (CS) The process $\{L(t), \ t\geq 0\}$ is stochastically continuous at zero. \\

\Bin \textbf{NB.} Each marginal law $K(t)$ is an infinitely decomposable (\textit{idecomp}) law.\\

\Bin \textbf{Examples}. \label{classical-examples} Here are three examples. (1) The Brownian process $(B_t)_{t\geq 0}$ with $K(t)=\mathcal{N}(0,t)$, $t\geq 0$; (2) The simple Poisson process $(N_t)_{t\geq 0}$ with $K(t)=\mathcal{P}(t)$, $t\geq 0$; (3) The $\gamma(a,b)$-process $(G_{a,b,t})_{t\geq 0}$ with  $K(t)=\gamma(at,b)$.\\

\Ni For L\'evy processes, we mainly refer to \cite{applebaum}.\\

\Bin \subsection{Motivation and aim} \label{ss_01_04}$ $\\

\Ni After this round-up on the key objects of that paper, we precise that we aim at doing some random analysis of the L\'evy process $\{Y^{\ast}(t), t \geq 0\}$ whose margins are defined by

$$
Y(t) \sim K(t)=\textit{gamma-gh}(at,b,0,\mu t, \sigma).
$$

\Bin The existence of such L\'evy processes is justified in \cite{niang-lo22}. Because of the more involvement of $GHD$ laws in Finance, we see it appropriate to have a detailed analysis of L\'evy processes  using \textit{GHD} margins should be done at the level of the details of Brownian motions for example. There are more studies on these stochastic processes in general. The specific studies we are going to do will help as example and counter-examples in those investigations.\\

\Ni We recall that we already have applied the Kolmogorov Existence Theorem (\textit{KET}) to prove the existence of L\'evy processes.\\

\Ni Hereafter, we begin with the \textit{gamma-gh}$(a,b,0,\mu , \sigma)$ L\'evy process because of the explicit form of the characteristic function.\\

\Ni (a) The non-where Totally Bounded variation of $T$.\\

\Ni (a) The non convergence of the Quadratic variation on compact sets in $L^2$.\\

\Ni (c) The empirical construction of $Y$ providing graphical illustrations.\\

\section{Study of the Variation of \textit{gamma-gh} Le\'vy processes}

\Ni Let $\mathcal{P}ar[a,b]$ be the collection of all finite partitions of $[0,T]$, $T>0$, of the form $\pi=(0=t_0<t_1< \cdots < t_k=T)$, $k\geq 1$. We define the variation of $Y^{\ast}$ on $\pi$ by 

$$
V_T(\pi)=\sum_{j=0}^{k-1} \left|Y^{\ast}(t_{j+1})-Y^{\ast}(t_j)\right|.
$$

\Bin The total variation of $Y^{\ast}$ on $[0,T]$ is defined as

$$
V_T=\sup_{\pi\in \mathcal{P}ar[a,b]} V_T(\pi).
$$

\Bin The stochastic process $Y$ is of bounded variation on some measurable subset $\Omega_0\subset \Omega$ if for any $\omega \in \Omega_0$\, $V_T(\omega)<+\infty$. If $\mathbb{P}(\Omega_0)>0$, a path-wise Riemann-Stieltjes integration theory of continuous stochastic processes $f(t,\omega)$, \textit{i.e.} treating integrals of the form,

$$
\Omega_0 \in \omega \rightarrow \int_{0}^{T} f(t,\omega) \ dY(t,\omega),
$$ 

\Bin is meaningful and interesting. Otherwise, \textit{i.e.} if $Y$ is non-where of bounded variation, that is 

$$
\mathbb{P}(V_T<\infty)=0,
$$

 then a Riemann-Stieltjes integration theory leads to a dead end. Methods of stochastic integration should be tried instead. The last described situation actually happens.\\

\Bin This is the commonly formulation of the problem. However, a firm look of the proof of that result (as given for example in \cite{mestuto-ang}) allows to see that we need less restrictive condition than the boundedness of total variation. We need only what we call the asymptotic boundedness of total variation we denote below. The full details related to that notion is to be found in \cite{ABV2022}.

\Ni Let $k\geq 1$ and  $\pi_{k}=(0=t_0^{(k)}<t_1^{(k)}< \cdots < t_{\ell(k)}^{(k)}=T)$, be a any finite partition of  $[0,T]$, \ $T>0$, with

$$
\ell(k)\rightarrow +\infty, \ \ m(\pi_{k})=\max_{0\leq j \leq \ell(k)-1} \left(t_{j+1}^{(k)}-t_j^{(k)}\right) \rightarrow 0 \ as \ k\rightarrow +\infty.  
$$

\Bin Let $F: [a,\ b]\rightarrow \mathbb{R}$ be a mapping. The variation of $F$ on $\pi_{k}$ is

$$
V_k = \sum_{j=0}^{\ell(k)-1} \left|F(t_{j+1}^{(k)})-F(t_j^{(k)})\right|.
$$ 

\Bin We say that $F$ has an  if and only if there exists a finite real number $M>0$ such that

$$
\limsup_{n\rightarrow +\infty} V_k<M
$$

\Bin whatever be the sequence of partitions $(V_k)_{k\geq 1}$ as $ell(k)\rightarrow +\infty$ and $ m(\pi_{k})asymptotic bounded variation$. It is shown in \cite{ABV2022} the asymptotic boundedness of variation is enough for making integrable all real-valued continuous functions on $[a, \ ]$ in the sense of the Riemann-Stieltjes approach.
\Bin We have

\begin{theorem} \label{nafy_01} 
The \textit{gamma-gh}$(a,b,0,\mu , \sigma)$ L\'evy process is almost-surely of asymptotic bounded variation.
\end{theorem}

\Bin Instead of directly analyzing $Y^{\ast}$, we choose to study the centered process

$$
\{Y(t), \ \ t\geq 0\}=\{Y^{\ast}(t)-\mu t, \ \  t\geq 0\}.
$$

\Bin \textbf{Proof of Theorem \ref{nafy_01}}. It is clear that $Y$ and $Y^{\ast}$ are of bounded variation at the same time or none of them is. From \eqref{ghd-rep}, we have that:

\begin{eqnarray*}
\forall \ 0<s<t, \ \ Z(s,t)=Y(t)-Y(s) &\sim& \textit{gamma-gh}(a(t-s),b,0,0, \sigma)\\
&\sim& \sigma Z_{(s,t)} \ L_{(s,t)},
\end{eqnarray*}

\Bin where  $Z_{(s,t)} \sim \mathcal{N}(0,1)$ is independent of $L_{(s,t)} \sim \gamma(a(t-s), b)^{1/2}$. Let $k\geq 1$ and  $\pi_{k}=(0=t_0^{(k)}<t_1^{(k)}< \cdots < t_{\ell(k)}^{(k)}=T)$, be a partition of  $[0,T]$, \ $T>0$, with

$$
\ell(k)\rightarrow +\infty, \ \ m(\pi_{k})=\max_{0\leq j \leq \ell(k)-1} \left(t_{j+1}^{(k)}-t_j^{(k)}\right) \rightarrow 0 \ as \ k\rightarrow +\infty.  
$$

\Bin The variation of $Y$ on $\pi_{k}$ is

$$
V_k = \sum_{j=0}^{\ell(k)-1} \left|Y(t_{j+1}^{(k)})-Y(t_j^{(k)})\right|.
$$

\Bin Let us show that $V_k$ has a finite variance. Let, for $0\leq j\leq \ell(k)-1$,

$$
\Delta_j = |Y(t_{j+1}^{(k)})-Y(t_{j}^{(k)})| \ \ and \ \ \delta_j=t_{j+1}^{(k)}-t_j^{(k)}.
$$

\Bin Each $\Delta_j$, $j \in \{0,\cdots,\ell(k)-1\}$ is of the form  $\Delta_j=\sigma |Z_j| L_j$, $Z_j\sim \mathcal{N}(0,1)$ independent of  $L_j \sim \gamma(a \delta_j,\beta)^{1/2}$. We have

\begin{eqnarray*}
\mathbb{E} \Delta_j&=&\sigma \sqrt{\frac{2}{\pi}} \int_{0}^{+\infty} x^{1/2} f_{\gamma(a \delta_j,\beta)}(x) \ dx\\
&=&\sigma \sqrt{\frac{2}{\pi}} \times \frac{\beta^{a \delta_j}}{\Gamma(a \delta_j)} \int_{0}^{+\infty} x^{a \delta_j+1/2-1} e^{-\beta x} \ dx\\
&=&\sigma \sqrt{\frac{2}{\pi}} \times \frac{\beta^{a \delta_j}}{\Gamma(a \delta_j)} \frac{\Gamma(a \delta_j+1/2)}{\beta^{a \delta_j+1/2}}\\
&=&\sigma \sqrt{\frac{2}{\pi\beta}} \times \frac{\Gamma(a \delta_j+1/2)}{\Gamma(a \delta_j)}.
\end{eqnarray*}

\Bin Let us denote

$$
h_j=\frac{\Gamma(a \delta_j+1/2)}{\Gamma(a \delta_j)}, \ 0\leq j \leq \ell(k)-1
$$

\Bin and

$$
\Gamma(a \delta_j+1/2)=\int_{0}^{1} x^{a \delta_j-1/2} \ e^{-x} \ dx + \int_{1}^{+\infty} x^{a \delta_j-1/2} \ e^{-x} \ dx=:b_{1,j}+b_{2,j}.
$$

\Bin Then, accordingly to: $e^{-1}\leq e^{-x}\leq 1$, we have

$$
\frac{2/e}{1+2a \delta_j} \leq b_{1,j} \leq \frac{2}{1+2a \delta_j},
$$

\Bin \textit{i.e.},

$$
(2/e) (1+\overline{o}(1)) \leq b_{1,j} \leq 2 (1+\overline{o}(1)),
$$

\Bin where $\overline{o}(1)$ is a small \textit{o} that is uniform in  $j$ over $a \delta_j\leq a m(\pi_k)$. As well,

$$
b_{2,j}=\int_{1}^{+\infty} e^{a \delta_j \log x} x^{-1/2} \ e^{-x} \ dx
$$ 

\Bin and

$$
I_1=\int_{1}^{+\infty}  x^{-1/2} \ e^{-x} \ dx \leq b_{2,j}\leq \int_{1}^{+\infty} e^{a m(\pi_k) \log x} x^{-1/2} e^{-x} \ dx.
$$ 

\Bin The function $e^{a m(\pi_k) \log x} x^{-1/2} e^{-x} $ is bounded by  $e^{\log x} x^{-1/2} e^{-x}= x^{1/2} e^{-x}  \in L^1([1,+\infty])$ and converges to $ x^{-1/2} e^{-x}$ as  $k\rightarrow +\infty$ (with $m(\pi_k)\rightarrow 0$). Then

$$
b_{2,j} \rightarrow I_1, 
$$

\Bin \ uniformly in $j$. We then have

$$
((2/e)+I_1) (1+\overline{o}(1)) \leq \Gamma(a \delta_j+1/2) \leq (2+I_1) (1+\overline{o}(1)). 
$$

\Bin Similarly, we have  

$$
I_2=\int_{1}^{\infty} x^{-1} \ e^{-x} \ dx,
$$

\Bin and for each $j$,

$$
\frac{1/e}{a \delta_j}+I_2 (1+\overline{o}(1)) \leq \Gamma(a \delta_j) \leq \frac{1}{a \delta_j}+I_2 (1+\overline{o}(1)) 
$$

\Bin \textit{i.e.},

$$
\frac{1/e+I_2 a \delta_j (1+\overline{o}(1))}{a \delta_j}  \leq \Gamma(a \delta_j) \leq \frac{1+I_2 a \delta_j (1+\overline{o}(1))}{a \delta_j}.
$$

\Bin Finally, we have

\begin{equation*}
\frac{\left( 2/e+I_{1}\right) a \delta _{j}\left( 1+\overline{o}\left(
1\right) \right) }{\left( 1+I_{2}a \delta _{j}\right) \left( 1+%
\overline{o}\left( 1\right) \right) }\leq h_{j}\leq \frac{\left(
2+I_{1}\right) a \delta _{j}\left( 1+\overline{o}\left( 1\right)
\right) }{\left( 1/e+I_{2}a \delta _{j}\right) \left( 1+\overline{o}%
\left( 1\right) \right) },
\end{equation*}

\Bin which implies

\begin{equation*}
\frac{\left( 2/e+I_{1}\right) a \delta _{j}\left( 1+\overline{o}\left(
1\right) \right) }{\left( 1+I_{2}a \delta _{j}\right) }\leq h_{j}\leq 
\frac{\left( 2+I_{1}\right) a \delta _{j}\left( 1+\overline{o}\left(
1\right) \right) }{\left( 1/e+I_{2}a \delta _{j}\right) }.
\end{equation*}

\Bin Then we have

\begin{equation*}
\sigma \sqrt{\frac{2}{\pi \beta }}\frac{\left( 2/e+I_{1}\right) a
\delta _{j}\left( 1+\overline{o}\left( 1\right) \right) }{\left(
1+I_{2}a \delta _{j}\right) }\leq \mathbb{E}\Delta _{j}\leq \sigma 
\sqrt{\frac{2}{\pi \beta }}\frac{\left( 2+I_{1}\right) a \delta
_{j}\left( 1+\overline{o}\left( 1\right) \right) }{\left( 1/e+I_{2}a
\delta _{j}\right) }
\end{equation*}

\Bin Let $\underline{m}\left( \pi _{k}\right) =\min_{1\leq j \leq \ell(k)-1} \delta _{j}\leq m(\pi_k$, then

\begin{equation*}
\sigma \sqrt{\frac{2}{\pi \beta }}\frac{\left( 2/e+I_{1}\right) a
\delta _{j}\left( 1+\overline{o}\left( 1\right) \right) }{\left(
1+I_{2}a m\left( \pi _{k}\right) \right) }\leq \mathbb{E}\Delta
_{j}\leq \sigma \sqrt{\frac{2}{\pi \beta }}\frac{\left( 2+I_{1}\right)
a \delta _{j}\left( 1+\overline{o}\left( 1\right) \right) }{\left(
1/e+I_{2}a \underline{m}\left( \pi _{k}\right) \right) }
\end{equation*}

\Bin and

\begin{equation*}
\sigma \sqrt{\frac{2}{\pi \beta }}\frac{\left( 2/e+I_{1}\right) a
T\left( 1+\overline{o}\left( 1\right) \right) }{\left( 1+I_{2}a m\left(
\pi _{k}\right) \right) }\leq \mathbb{E}V_{k}\leq \sigma \sqrt{\frac{2}{\pi
\beta }}\frac{\left( 2+I_{1}\right) a T\left( 1+\overline{o}\left(
1\right) \right) }{\left( 1/e+I_{2}a \underline{m}\left( \pi
_{k}\right) \right) }.
\end{equation*}

\Bin Put $E_{1}=2/e+I_{1}$ and $E_{2}=\frac{2+I_{1}}{1/e}$. So, as $m\left( \pi _{k}\right) \rightarrow 0$, we get 

\begin{equation*}
\sigma \sqrt{\frac{2}{\pi \beta }}E_{1}a T\leq \mathbb{E}V_{k}\leq
\sigma \sqrt{\frac{2}{\pi \beta }}E_{2}a T
\end{equation*}

\Bin Now, let us compute $\mathbb{V}ar(V_k)$. Since the $\Delta_j$'s are independent, we have

$$
\mathbb{V}ar(V_k)=\sum_{j=0}^{\ell(k)-1} \mathbb{V}ar(\Delta_j). 
$$

\Bin But we have for each $j$,

$$
\mathbb{E} \Delta_j^2=\mathbb{E}(\sigma^2 Z_j^2 L_j^2)=\sigma^2 \frac{a \delta_j}{\beta}.
$$

\Bin Hence

\begin{eqnarray*}
&& \sigma^2 \frac{a \delta_j}{\beta}- E_2^2  \frac{2\sigma^2}{\pi\beta} (a\delta_j)^2 \ (1+\overline{o}(1))\\
&& \leq \mathbb{V}ar(\Delta_j)\\
&& \leq \sigma^2 \frac{a \delta_j}{\beta}- E_1^2  \frac{2\sigma^2}{\pi\beta}  (a\delta_j)^2 \ (1+\overline{o}(1)).
\end{eqnarray*}

\Bin We get, as $k\rightarrow +\infty$, 

$$
\sum_{j=0}^{\ell(k)-1} \delta_j^2\leq m(\pi_k) T \rightarrow 0.
$$

\Bin By summing over $j$, we arrive at

\begin{eqnarray*}
\sigma^2 \frac{a T}{\beta}- \overline{o}(1)\leq \mathbb{V}ar(V_k) \leq \sigma^2 \frac{a T}{\beta}- \overline{o}(1).
\end{eqnarray*}

\Bin Then

$$
\mathbb{V}ar(V_k)=\sigma^2 \frac{a T}{\beta} + \overline{o}(1).
$$

\Bin We conclude as follows. Let $K_1=\sigma \sqrt{\frac{2}{\pi \beta }}E_{2}a T$ and for $K_2=\sigma^2 \frac{a T}{\beta}$. Then for any $p\geq 1$

\begin{equation}
(\sup_{k\geq 1} V_k=+\infty)=\bigcap_{n\geq 1} \biggr(\sup_{k\geq 1} V_k \geq n\biggr)\subset \biggr(\sup_{k\geq 1} V_k \geq p\biggr)=A_p. \label{approx1}
\end{equation}

\Bin Next we have

$$
A_{p,q}=\biggr(\sup_{1\leq k \leq q} V_k \geq p\biggr) \nearrow A_p \ as \ q\nearrow +\infty.
$$

\Bin Now, for the we may superpose the partitions $(\pi_k)_{1\leq k\leq q}$ into the partition 

$$
\pi^\ast_{q}=(0=c^\ast_{1,q}<c^\ast_{2,q}<\cdots<c^\ast_{1,\ell^\ast(q)})
$$

\Bin Let us denote $V^\ast_{q}$ as the variation of $Y$ over the partition $\pi^\ast_{q}$. We can directly see that:\\

\Ni (a) For any $k \in [1,q]$, $m(\pi^\ast_{q})\leq m(\pi_k)$;\\

\Ni (b) For any $k \in [1,q]$, $V_k \leq V^\ast_{q}$.\\

\Ni Then we have

$$
\mathbb{P}(A_{p,q})\leq \mathbb{P}(V^\ast_{q}\geq p)\leq \frac{\mathbb{E}V^\ast_{q}}{p}\leq \frac{K_1}{p}.
$$

\Bin We let $q\nearrow +\infty$ to get

$$
\mathbb{P}(A_{p})\leq \mathbb{P}(V^\ast_{q}\geq p)\leq \frac{K_1}{p}.
$$

\Bin Next, we combine this with Inequality \eqref{approx1} and let $p\nearrow +\infty$ to have
\Bin Then

$$
\mathbb{P}(\sup_{k\geq 1} V_k=+\infty)=0.
$$

\Bin We conclude that

$$
\mathbb{P}(\limsup_{n\rightarrow +\infty} V_k)=0. \ \ \ \ \blacksquare
$$

\Bin Now, we move the quadratic variation.\\

\section{Quadratic Variation}

\Ni The Quadratic Variation of $Y$ over $[0, \ T]$, by definition, is the limit (in some sense), as $k\rightarrow +\infty$, of 

$$
V_Q(T,\pi_k)=\sum_{j=0}^{\ell(k)-1} (Y(t_{j+1}^{(k)})-Y(t_{j}^{(k)}))^2. 
$$

\Bin Here is the result.

\begin{theorem} \label{nafy_02} 
The Quadratic Variation $V_Q(T,\pi_k)$ does not converge to a constant in $L^2$  as $m(\pi_k)\rightarrow 0$. 
\end{theorem}

\Bin \textbf{Proof of Theorem \ref{nafy_02}}. We have

$$
V_Q(T,\pi_k)=\sum_{j=0}^{\ell(k)-1} \Delta_j^2 
$$

\Bin and

$$
\mathbb{E} V_Q(T,\pi_k) = \sum_{j=0}^{\ell(k)-1} \mathbb{E}\Delta_j^2=\sum_{j=0}^{\ell(k)-1} \sigma^2 \mathbb{E}Z_j^2 \mathbb{E}\left(\gamma(a \delta_j,\beta)\right).
$$

\Bin Hence

$$
\mathbb{E} V_Q(T,\pi_k)=\sum_{j=0}^{\ell(k)-1} \sigma^2 a \delta_j/\beta=a \sigma^2 T/\beta. 
$$

\Bin By independent of increments,

$$
\mathbb{V}ar\left(V_Q(T,\pi_k)\right)=\sum_{j=0}^{\ell(k)-1} \mathbb{V}ar(\Delta_j^2)
$$

\Bin and we already know that

$$
\mathbb{E} \Delta_j^2=\sigma^2 a \delta_j/\beta.
$$

\Bin Now, we have

$$
\mathbb{E} \Delta_j^4=\mathbb{E} \biggr(\sigma Z_j \gamma(a \delta_j,\beta)^{1/2}\biggr)^4= 3\sigma^4  \frac{(a \delta_j)(a \delta_j+1)}{\beta^2},
$$

\Bin \textit{i.e.},

$$
\mathbb{E} \Delta_j^4= 3\sigma^4 \frac{(a \delta_j)^2}{\beta^2}+3\sigma^4 \frac{(a \delta_j)}{\beta^2}.
$$

\Bin This leads to

$$
\mathbb{V}ar(\Delta_j^2)=2\sigma^4 \frac{(a \delta_j)^2}{\beta^2}+3\sigma^4 \frac{(a \delta_j)}{\beta^2}.
$$

\Bin We conclude that

$$
\mathbb{V}ar(V_Q(T,\pi_k))=2 \frac{\sigma^4 a^2 }{\beta^2}\sum_{j=0}^{\ell(k)-1} \delta_j^2 + 3 \frac{a\sigma^4 T}{\beta^2}. 
$$

\Bin Since

$$
0\leq \sum_{j=0}^{\ell(k)-1} \delta_j^2\leq m(\pi_k) T \rightarrow 0,
$$

\Bin we get 

$$
\mathbb{V}ar(V_Q(T,\pi_k))\rightarrow C= 3 \frac{a\sigma^4 T}{\beta^2}.
$$

\Bin By means of Huygens Formula, for any  $c\in \mathbb{R}$, for $\|\circ\|^2=\mathbb{E} (\circ^2)$, we obtained

$$
\|V_Q(T,\pi_k)-c\|^2= \|V_Q(T,\pi_k)-\mathbb{E} V_Q(T,\pi_k)\|^2 + (c-\mathbb{E} V_Q(T,\pi_k))^2\geq \biggr(3 \frac{a\sigma^4 T}{\beta^2} + o(1)\biggr).
$$

\Bin Hence $V_Q(T,\pi_k)$ does not converge in  $L^2$. $\blacksquare$\\

\Bin \textbf{A quick comparison with the Brownian process}. The \textit{gamma-gh} stochastic process from with the Brownian motion and diverges regarding the quadratic variation and the asymptotic boundedness of variation.\\

\Bin Now, let us move to the empirical construction that will allow us to get graphical representations.\\

\section{Empirical construction}

\Ni Let $T>0$. For each $n\geq 1$, let $X_{j,n}$, $1\leq j\leq n$, be independent random variables following the same distribution \textit{gamma}$(aT/n,b,0,\mu T/n,\sigma)$. By the infinite divisibility process, we easily get that

$$
\sum_{1\leq j \leq n} X_{j,n}=:Z_n \sim \textit{gamma}(a T,b,0,\mu T,\sigma) \textit{ for } n\geq 1.
$$

\Bin Now define, for $n\geq 1$, for any $T>0$

\begin{equation}
Y^{\ast}_n(t)= \biggr( \ \sum_{1\leq j \leq [nt/T]} X_{j,n}\biggr) 1_{([nt/T]\geq 1)}, \ \ t\in [0,T]. \label{levyEmp}
\end{equation}

\Bin In the frame of weak convergence of stochastic processes in $\ell^{\infty}(0,T)$, we have to check that the finite distribution convergence of $Y^{\ast}_n$. For $k\geq 1$, for $0=t_0<t_1<\ldots<t_k=T$, we have that the vectors, $1\leq j\leq k$,

$$
\Delta_j=Y^{\ast}_n(t_j)-(t_{j-1})= \sum_{[nt_{j-1}/T]+1}^{[nt_{j}/T]} X_{j,n},
$$

\Bin are independent. For each $j \in \{1,\ldots,k\}$, by using \eqref{gammagh-fchaf}, we clearly have

\begin{eqnarray*}
\psi_{\Delta_j}(u)&=&\exp\biggr(iu\mu T([nt_{j}/T]-[nt_{j-1}/T])/n\biggr) \biggr(1+ \frac{\sigma^2u^2}{2b}\biggr)^{-a T([nt_{j}/T]-[nt_{j-1}/T])/n}\\
&\rightarrow &\exp\biggr(i\mu (t_{j}-t_{j-1}) u\biggr) \biggr(1+ \frac{\sigma^2u^2}{2b}\biggr)^{-a(t_{j}-t_{j-1})},
\end{eqnarray*}

\Bin as $n\rightarrow +\infty$. So

$$
\Delta_j \rightsquigarrow  \textit{gamma}(a(t_{j}-t_{j-1}),b,0,\mu (t_{j}-t_{j-1}),\sigma), \textit{ as } n \rightarrow +\infty.
$$

\Bin By the Slutsky's theorem, we conclude that the finite distribution of $Y^{\ast}_n$ weakly converges those of \textit{gamma}$(at,b,0,\mu t,\sigma)$-L\'evy process. To formally have the full weak convergence, we have to show the asymptotic, which is to be done following the methods in \cite{niangLTB2022}. We will not do it here. Rather, we focus one approximation graphical representations using those of $Y^{\ast}_n$.\\

%gamma-levy1
%gamma-levy3
%gamma-levy10
%gamma-levy-brownian
%gamma-levy0p5

\Bin We propose approximations of the centered \textit{gamma}$gh(at,b,c,\mu t,\sigma)$-L\'evy processes by $Y_n$ in $[0, \ T]$. We take $n=500$, $b=\mu=1$, $\sigma=0.5$. We give four different values to $a$, \textit{i.e.} $a \in \{0.5, 1, 3, 10\}$. Finally, we represent the Brownian motion on $[0, \ T]$ as the weak limit of the sequence of stochastic processes

$$
Z^{n}(t)=\frac{1}{\sqrt{n}} \sum_{j=1}^{[nt/T]} Z_j,
$$

\begin{figure}
	\centering
		\includegraphics[width=0.75\textwidth]{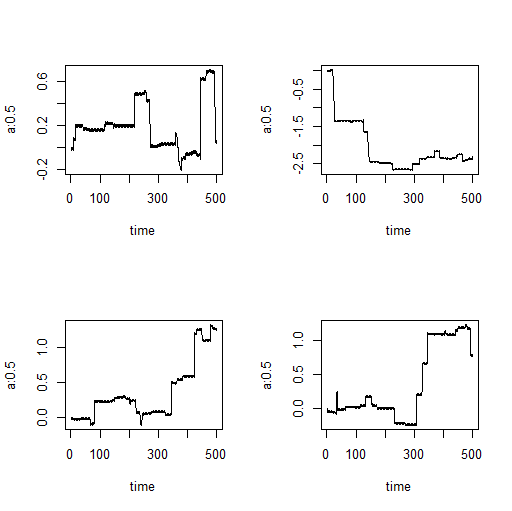}
	\caption{Construction of the gamma-gh L\'evy process for \\$a=0.5;b=1;\mu=1;\sigma=0.5$}
	\label{fig0p5}
\end{figure}

\begin{figure}
	\centering
		\includegraphics[width=0.75\textwidth]{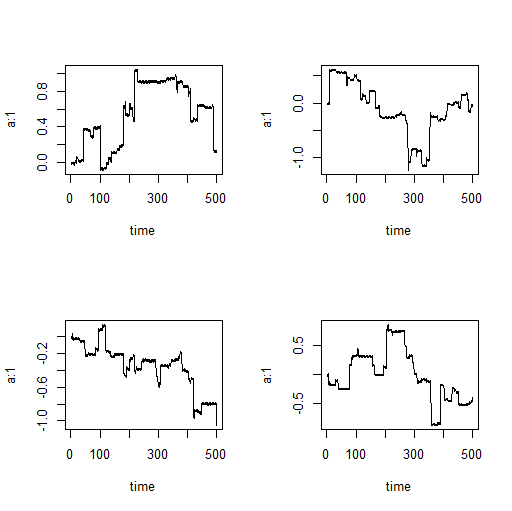}
	\caption{Construction of the gamma-gh L\'evy process for \\$a=1;b=1;\mu=1;\sigma=0.5$}
	\label{fig1}
\end{figure}

\begin{figure}
	\centering
		\includegraphics[width=0.75\textwidth]{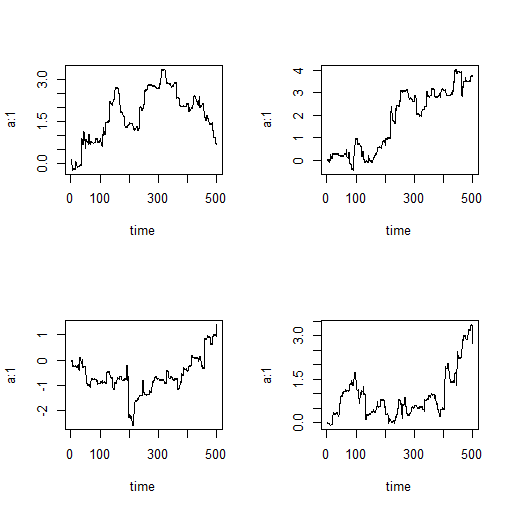}
	\caption{Construction of the gamma-gh L\'evy process for \\$a=3;b=1;\mu=1;\sigma=0.5$}
	\label{fig3}
\end{figure}

\begin{figure}
	\centering
		\includegraphics[width=0.75\textwidth]{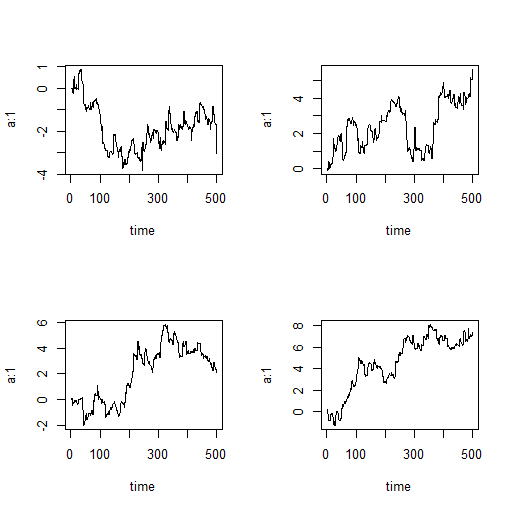}
	\caption{Construction of the gamma-gh L\'evy process for \\$a=10;b=1;\mu=1;\sigma=0.5$}
	\label{fig10}
\end{figure}

\begin{figure}
	\centering
		\includegraphics[width=0.75\textwidth]{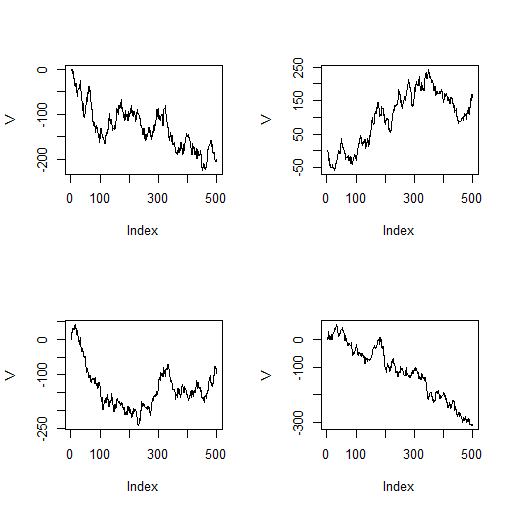}
	\caption{Construction of the Brownian motion}
	\label{figbrownian}
\end{figure}

\Bin where the $Z_i$ are independent $\mathcal{N}(0,T)$ random variables. In Figures \ref{fig0p5} and \ref{fig1} 
(pages \pageref{fig0p5} and \pageref{fig1}), for small values
$a \in \{0.5, 1\}$, the \textit{gamma}$gh(at,b,0,\mu t,\sigma)$-L\'evy processes have shapes very different that of the Brownian motion in Figure \ref{figbrownian} (page \pageref{figbrownian}). But as $a$ gets bigger, like $a \in \{3, 10\}$, the shapes of the \textit{gamma}$gh(at,b,c,\mu t,\sigma)$-L\'evy processes in Figures \ref{fig3} and \ref{fig10} (pages \pageref{fig3} and \pageref{fig10}) and the Brownian motion seems to be alike. It was important for us to see what might be like paths of the studied process.\\

\section{Conclusion} 

\Ni This study is the first on a series of random analysis studies on $gh(at,b,c,\mu t,\sigma)$-L\'evy processes. Different aspects of stochastic analysis will be investigated later.\\

\Bin \textbf{Acknowledgment}.

\newpage
\newpage

\end{document}